\documentclass{article}
\usepackage{amsfonts}
\usepackage{amsmath}
\usepackage{amssymb}
\usepackage{setspace}
\input{tcilatex}
\begin{document}
\title{The group of rational points on the Holm curve is torsion-free}
\author{Fredrick M. Nelson}
\maketitle
\tableofcontents
\begin{abstract}
Using the division polynomials for elliptic curves in Weierstrass form, it is shown that the group of rational points on the curve $H:k(y^{3}-y)=l(x^{3}-x)$ is torsion-free.
\end{abstract}
\doublespacing
\section{Introduction}
\paragraph{}The equation of the Holm curve is $H:k(y^{3}-y)=l(x^{3}-x)$ where $k$ and $l$
are distinct, positive, relatively prime, square-free integers. Holm \cite{holm}
used $H$ to show how to construct rational right triangles whose areas have ratio $k/l$.
\paragraph{} A positive integer $n$ is a
congruent number if there exists a rational right triangle with area $n$. \
Rajan and Ramarason \cite{rajan-ramaroson} used $H$ to show that there are
infinitely many congruent numbers with ratio $k/l$.
\newline\newline
\paragraph{} It will be shown that the group of rational points on $H$ has no points of finite order.
\section{Definitions}
\paragraph{} $\mathbb{Z}$ is the set of integers and $\mathbb{Q}$ is the set of rational numbers.
\newline
If $p$ is prime and $x\in \mathbb{Q}-\{0\}$, then $v_{p}\left(
x\right) $ is the index of $p$ in the prime factorization of $x$.
\newline
$H(\mathbb{Q})=\left\{\text{$(x,y)\in\mathbb{Q}^2:k(y^3-y)=l(x^3-x)$}
\right\}$ and
\newline
$E(\mathbb{Q})=\left\{\text{$(x,y)\in\mathbb{Q}
^2:y^2=x^{3}-3k^{2}l^{2}x+k^{2}l^{2}( k^{2}+l^{2}))$}\right\}\cup\left\{{
\infty}\right\}$.
\newline
If $P = (a,b)\in E(\mathbb{Q})$, then $x(P) = a$ and $y(P) = b$.
\section{Isomorphism}
\paragraph{} By Mordell's theorem, $H(\mathbb{Q})$ and $E(\mathbb{Q})$ are finitely
generated abelian groups. Choose $(0,0)$ as the identity element for $H(
\mathbb{Q})$ and the point at infinity $\infty$ as the identity element for $E(\mathbb{Q})$.
\paragraph{}Then the map $\gamma :H(\mathbb{Q})\rightarrow E(\mathbb{Q})$ given by 
\begin{equation}
\gamma (\allowbreak x,\allowbreak y)=\left( \frac{kl\left( kx-ly\right) }{
lx-ky},\frac{kl\left( k^{2}-l^{2}\right) }{lx-ky}\right) 
\end{equation}
is an isomorphism. The inverse map is given by 
\begin{equation}
\gamma ^{-1}(\allowbreak x,y)=\left( \frac{k\left( x-l^{2}\right) }{y},\frac{
l\left( x-k^{2}\right) }{y}\right) .
\end{equation}
\section{Three Lemmas}
\paragraph{} $E(\mathbb{Q})$ is in Weierstrass form so the points of finite order of $E(\mathbb{Q})$ are integral by the Nagell-Lutz Theorem.
Since $\gamma$ is an isomorphism, we need only show that each non-identity
element of $E(\mathbb{Q})$ has a multiple which is not integral.  For this, it suffices to prove the following three lemmas.\newline\newline
Lemma 1: \newline
If $\left( x,y\right) $ is an integral point in $E\left( 
\mathbb{Q}
\right) $ and $2$ divides $k$ or $l$, then
\newline\newline $v_{2}\left(
x\left( 2(\allowbreak x,\allowbreak y)\right) \right) =\left\{ 
\begin{array}{rrr}
0 & \text{if} & v_{2}\left( \allowbreak x\allowbreak \right) \geq 2,\\ 
2 & \text{if} & v_{2}\left( \allowbreak x\allowbreak \right) =1,\\ 
 -2\text{ } & \text{if } & v_{2}\left( \allowbreak x\allowbreak \right) =0.
\end{array}
\right. $\newline\newline
Hence, $8(\allowbreak x,\allowbreak y)$ is not an integral point of $E(\mathbb{Q})$.
\newline\newline
Lemma 2: \newline
If $\left( x,y\right) $ is an integral point in $E\left( 
\mathbb{Q}
\right) $ and $p$ is an odd prime dividing $k$ or $l$, then
\newline\newline
$v_{p}\left( x\left( 3\left( x,y\right) \right) \right) \leq \left\{ 
\begin{array}{lll}
-2\text{ } & \text{if } & p=3\text{,} \\ 
0 & \text{if } & \text{otherwise.}
\end{array}
\right. $
\newline\newline
Hence, $3(\allowbreak x,\allowbreak y)$ is not an integral point of $E(\mathbb{Q})$.
\newline\newline
Lemma 3: $\ $If $\left( x,y\right) $ is an integral point in $E\left( 
\mathbb{Q}
\right) $ and $p\geq 5$ is a prime dividing $k$ or $l$, then $v_{p}\left(
x\left( 3p\left( x,y\right) \right) \right) \leq -2$.
\newline
Hence, $3p(\allowbreak x,\allowbreak y)$ is not an integral point of $E(\mathbb{Q})$.
\newline\newline
\section{Division Polynomials}
\paragraph{} Let $a,b\in \mathbb{Z}$ where $4a^{3}+27b^{2}\neq 0$.
Then $y^{2}=x^{3}+ax+b$ is an elliptic curve in Weierstrass form.
\newline\newline
The division polynomials $\psi _{n}$ are defined as follows:
\newline
$\psi _{0}=0$,
\newline
$\psi _{1}=1$,
\newline
$\psi _{2}=2y$,
\newline
$\psi _{3}={3}x^{4}+6ax^{2}+12bx-a^{2}$,
\newline
$\psi _{4}=y\left(
4x^{6}+20ax^{4}+80bx^{3}-20a^{2}x^{2}-16abx-4a^{3}-32b^{2}\right) $,
\newline
$ \psi _{2n+1}=\psi _{n+2}\psi _{n}^{3}-\psi _{n-1}\psi
_{n+1}^{3}$ for $n\geq 2$, and 
\newline
$\psi _{2n}=\frac{\psi _{n}\left( \psi _{n+2}\psi
_{n-1}^{2}-\psi _{n-2}\psi _{n+1}^{2}\right) }{2y}$ for $n\geq 3$.
\newline\newline
Auxiliary functions $\phi _{n}$ and $\omega _{n}$ are defined from $\psi _{n}$ by
\newline
$\phi _{n}=x\psi _{n}^{2}-\psi _{n-1}\psi _{n+1}$ for $n\geq 1$ and
\newline
$\omega _{n}=\frac{\psi _{n-1}^{2}\psi _{n+2}-\psi
_{n-2}\psi _{n+1}^{2}}{4y}$ for $n\geq 2$.
\newline\newline
From these definitions it follows that $n\left( x,y\right) =\left( \frac{\phi _{n}}{\psi _{n}^{2}
},\frac{\omega _{n}}{\psi _{n}^{3}}\right) $ if $n\geq 2$ and $\left( x,y\right)$ is a point on the curve  $y^{2}=x^{3}+ax+b$.
\newline\newline
The following proposition follows from the properties of division polynomials described in \cite{mckee}.
\newline\newline
Proposition
\newline
(i) $2\left( x,y\right) =\left( \frac{\allowbreak x^{4}-2ax^{2}-8bx+a^{2}}{
\allowbreak 4y^{2}},\frac{x^{6}+5ax^{4}+20bx^{3}-5a^{2}x^{2}-4ab\allowbreak
x-a^{3}-8b^{2}}{\allowbreak 8y^{3}}\right) $.
\newline
(ii) $x\left( 3\left( x,y\right) \right) =\frac{\allowbreak
x^{9}-12ax^{7}-96bx^{6}+30a^{2}x^{5}-24abx^{4}+\allowbreak 12\left(
3a^{3}+4b^{2}\right) x^{3}+48a^{2}b\allowbreak x^{2}+\allowbreak 3a\left(
3a^{3}+32b^{2}\right) x+8b\left( a^{3}+8b^{2}\right) }{\left( {3}
x^{4}+6ax^{2}+12bx-a^{2}\right) ^{2}}$.
\newline
(iii) The degree and leading coefficient of ${\psi _{n}^{2}}$ are $n^{2}-1$
and $n^{2}$, respectively.
\newline
(iv) The degree and leading coefficient of $\phi _{n}$ are $n^{2}$ and $1$,
respectively.
\newline
(v) The coefficient of $x^{n^{2}-2}$ in ${\psi _{n}^{2}}$ and the
coefficient of $x^{n^{2}-1}$ in $\phi _{n}$ are both $0$.
\newline
(vi) If a positive integer $d$ divides both $a$ and $b$, then $d$ divides
both $\psi _{n}^{2}-n^{2}x^{n^{2}-1}$ and $\phi _{n}-x^{n^{2}}$.
\newline
(vii) Let $d,n,w,z$ be integers with $d\geq 1,n\geq 2$ and $w,z\neq 0$. \
Let $x=\frac{w}{z}$.
\newline
Then $z^{n^{2}-3}\left( \psi _{n}^{2}-n^{2}x^{n^{2}-1}\right) $ and $
z^{n^{2}-2}\left( \phi _{n}-x^{n^{2}}\right) $ are integers.  Furthermore, if $d$ divides both $a$ and $b$ then $d$ divides both $
z^{n^{2}-3}\left( \psi _{n}^{2}-n^{2}x^{n^{2}-1}\right) $ and $
z^{n^{2}-2}\left( \phi _{n}-x^{n^{2}}\right) $.
\section{Proofs}
For the remainder of the paper, let $a=-3k^{2}l^{2}$ and $b=k^{2}l^{2}\left( k^{2}+l^{2}\right)$.   
\newline
Then $E\left(\mathbb{Q}\right) =\left\{ \left( x,y\right) \in 
\mathbb{Q}
^{2}:y^{2}=x^{3}+ax+b\right\} \cup \left\{ \infty \right\} $.
\newline
Note that $k^{2}l^{2}$ is a positive integer that divides both $a$ and $b$.
\newline\newline
Also, if $\left( x,y\right)\in E\left(\mathbb{Q}\right)-\left\{{
\infty}\right\}$, then both $x$ and $y$ are nonzero.  For if $x=0$, then $y^{2}=k^{2}l^{2}\left( k^{2}+l^{2}\right) $, $k^{2}+l^{2}$
is a square, and one of $k$ and $l$ has a square factor, a contradiction.  Hence, $x\neq 0$.
If $y=0$, then $\exists \left( u,v\right) \in E\left( 
\mathbb{Q}
\right) $ such that $\left( \frac{kl\left( k\allowbreak u-l\allowbreak
v\right) }{l\allowbreak u-k\allowbreak v},\frac{kl\left( k^{2}-l^{2}\right) 
}{l\allowbreak u-k\allowbreak v}\right) =\gamma (\allowbreak \allowbreak
u,\allowbreak \allowbreak v)=\left( x,0\right) $, $\frac{kl\left(
	k^{2}-l^{2}\right) }{l\allowbreak u-k\allowbreak v}=0$, $kl\left(
k^{2}-l^{2}\right) =0$, $k^{2}-l^{2}=0$, and $k=l$, a contradiction$.$
Hence, $y\neq 0$.
\newline\newline
Proof of Lemma 1
\newline
We have $k,l\in 
\mathbb{Z}
$, $k,l>0$, $k\neq l$, $\gcd \left( k,l\right) =1$, and $k$ and $l$ are
square-free.
We may assume without loss of generality that $2|k$.
Then $v_{2}\left( k\right) =1$ and $v_{2}\left( l\right) =0$.
Since $v_{2}\left( k\right) =1$, $\exists j\in \mathbb{Z}
-\left\{ 0\right\} $ such that $k=2j$ and $v_{2}\left( j\right) =0$.
Then $y^{2}=x^{3}-3k^{2}l^{2}x+k^{2}l^{2}\left( k^{2}+l^{2}\right)
=x^{3}-3\left( 2j\right) ^{2}l^{2}x+\left( 2j\right) ^{2}l^{2}\left( \left(
2j\right) ^{2}+l^{2}\right) =\allowbreak
x^{3}-12j^{2}l^{2}x+4j^{2}l^{2}\left( 4j^{2}+l^{2}\right) $ and
$x\left( 2\left( x,y\right) \right) =\allowbreak \frac{\allowbreak
\allowbreak x^{4}+6k^{2}l^{2}\allowbreak x^{2}\allowbreak -8k^{2}l^{2}\left(
k^{2}+l^{2}\right) \allowbreak x+9k^{4}l^{4}}{4\allowbreak y^{2}}
=\allowbreak \frac{\allowbreak \allowbreak x^{4}+6\left( 2j\right)
^{2}l^{2}\allowbreak x^{2}\allowbreak -8\left( 2j\right) ^{2}l^{2}\left(
\left( 2j\right) ^{2}+l^{2}\right) \allowbreak x+9\left( 2j\right) ^{4}l^{4}
}{4\left( \allowbreak x^{3}-12j^{2}l^{2}x+4j^{2}l^{2}\left(
4j^{2}+l^{2}\right) \right) }=\frac{\allowbreak
x^{4}+24j^{2}l^{2}x^{2}-32j^{2}l^{2}\left( 4j^{2}+l^{2}\right)
x+144j^{4}l^{4}}{4\allowbreak \left( \allowbreak
x^{3}-12j^{2}l^{2}x+\allowbreak 4j^{2}l^{2}\left( 4j^{2}+l^{2}\right)
\right) }$.
\newline
Case 1: \ $v_{2}\left( \allowbreak x\allowbreak
\right) \geq 2$.
Then $\left( \exists n\geq 2\right) \left( \exists u\in 
\mathbb{Z}
-\left\{ 0\right\} \right) $ such that $x=2^{n}u$ and $v_{2}\left(
\allowbreak u\allowbreak \right) =0$.
Then $x\left( 2\left( x,y\right) \right) =\frac{\allowbreak \left(
2^{n}u\right) ^{4}+24j^{2}l^{2}\left( 2^{n}u\right) ^{2}-32j^{2}l^{2}\left(
4j^{2}+l^{2}\right) \left( 2^{n}u\right) +144j^{4}l^{4}}{4\allowbreak \left(
\allowbreak \left( 2^{n}u\right) ^{3}-12j^{2}l^{2}\left( 2^{n}u\right)
+\allowbreak 4j^{2}l^{2}\left( 4j^{2}+l^{2}\right) \right) }=\frac{
\allowbreak 2^{n+1}\left( 2^{3n-5}u^{4}+2^{\allowbreak
n-2}3j^{2}l^{2}u^{2}-j^{2}l^{2}\left( 4j^{2}+l^{2}\right) \allowbreak
u\right) +9j^{4}l^{4}}{\allowbreak 2^{n}\left(
2^{2n-2}u^{3}-3j^{2}l^{2}u\right) +j^{2}l^{2}\left( 4j^{2}+l^{2}\right) }$,
\newline
the numerator and denominator of $x\left( 2\left( x,y\right) \right) $ are both odd,
and $v_{2}\left( x\left( 2\left( x,y\right) \right) \right) =0$.
\newline
Case 2: \ $v_{2}\left( \allowbreak x\allowbreak \right) =1$.
Then $\exists u\in 
\mathbb{Z}
-\left\{ 0\right\} $ such that $x=2u$ and $v_{2}\left( \allowbreak
u\allowbreak \right) =0$.
\newline
Since $v_{2}\left( \allowbreak u\allowbreak \right) =v_{2}\left( j\right)
=v_{2}\left( l\right) =0$, $\exists $ $f,g,h\in 
\mathbb{Z}
$ such that $u=2f+1$, $j=2g+1$, and $l=2h+1$.
\newline
Then $x\left( 2\left( x,y\right) \right) =\allowbreak \frac{\allowbreak
\allowbreak \left( 2u\right) ^{4}+24j^{2}l^{2}\allowbreak \left( 2u\right)
^{2}-32j^{2}l^{2}\left( 4j^{2}+l^{2}\right) \allowbreak \left( 2u\right)
+144j^{4}l^{4}}{4\allowbreak \left( \allowbreak \allowbreak \left( 2u\right)
^{3}-12j^{2}l^{2}\allowbreak \left( 2u\right) +\allowbreak 4j^{2}l^{2}\left(
4j^{2}+l^{2}\right) \right) }$
\newline
$=\allowbreak \frac{\allowbreak \allowbreak \left( 2\allowbreak \left(
2f+1\right) \right) ^{4}+24\allowbreak \left( 2g+1\right) ^{2}\allowbreak
\left( 2h+1\right) ^{2}\allowbreak \left( 2\allowbreak \left( 2f+1\right)
\right) ^{2}-32\allowbreak \left( 2g+1\right) ^{2}\allowbreak \left(
2h+1\right) ^{2}\left( 4\allowbreak \left( 2g+1\right) ^{2}+\allowbreak
\left( 2h+1\right) ^{2}\right) \allowbreak \left( 2\allowbreak \left(
2f+1\right) \right) +144\allowbreak \left( 2g+1\right) ^{4}\allowbreak
\left( 2h+1\right) ^{4}}{4\allowbreak \left( \allowbreak \allowbreak \left(
2\allowbreak \left( 2f+1\right) \right) ^{3}-12\allowbreak \left(
2g+1\right) ^{2}\allowbreak \left( 2h+1\right) ^{2}\allowbreak \left(
2\allowbreak \left( 2f+1\right) \right) +\allowbreak 4\allowbreak \left(
2g+1\right) ^{2}\allowbreak \left( 2h+1\right) ^{2}\left( 4\allowbreak
\left( 2g+1\right) ^{2}+\allowbreak \left( 2h+1\right) ^{2}\right) \right) }$.
\newline
This expression for $x\left( 2\left( x,y\right) \right)$ simplifies to a fraction whose numerator is $4$ times an odd number and whose denominator is odd.  Hence, $v_{2}\left( x\left( 2\left( x,y\right) \right) \right) =v_{2}\left(
4\right) =2$.
\newline
Case 3: \ $v_{2}\left( \allowbreak x\allowbreak \right) =0$.
\newline
Then $\exists u\in \mathbb{Z}$ such that $x=2u+1$ and $x\left( 2\left( x,y\right) \right) =\frac{
\allowbreak x^{4}+24j^{2}l^{2}x^{2}-32j^{2}l^{2}\left( 4j^{2}+l^{2}\right)
x+144j^{4}l^{4}}{4\allowbreak \left( \allowbreak
x^{3}-12j^{2}l^{2}x+\allowbreak 4j^{2}l^{2}\left( 4j^{2}+l^{2}\right)
\right) }$
\newline
$=\frac{\allowbreak \left( 2u+1\right) ^{4}+24j^{2}l^{2}\left( 2u+1\right)
^{2}-32j^{2}l^{2}\left( 4j^{2}+l^{2}\right) \left( 2u+1\right) +144j^{4}l^{4}
}{4\allowbreak \left( \allowbreak \left( 2u+1\right) ^{3}-12j^{2}l^{2}\left(
2u+1\right) +\allowbreak 4j^{2}l^{2}\left( 4j^{2}+l^{2}\right) \right) }$
\newline
$=\allowbreak \frac{8\left(
u+3j^{2}l^{2}-4j^{2}l^{4}-16j^{4}l^{2}+18j^{4}l^{4}+3u^{2}+4u^{3}+2u^{4}+12j^{2}l^{2}u^{2}+12j^{2}l^{2}u-8j^{2}l^{4}u-32j^{4}l^{2}u\right) \allowbreak +1
}{4\left( \allowbreak 2\left(
3u-6j^{2}l^{2}+2j^{2}l^{4}+8j^{4}l^{2}+6u^{2}+4u^{3}-12j^{2}l^{2}u\right)
+1\right) }\allowbreak $,
\newline
the numerator of $x\left( 2\left( x,y\right) \right) $ is odd, the denominator of $x\left( 2\left( x,y\right) \right) $ is $4$ times an odd number, and $v_{2}\left( x\left( 2\left( x,y\right) \right) \right) =-2$.
\newline
If $v_{2}\left( \allowbreak x\allowbreak \right) =0$, then $v_{2}\left(
x\left( 2_{E}(\allowbreak x,\allowbreak y)\right) \right) =-2$ and $
2_{E}(\allowbreak x,\allowbreak y)$ is not an integral point.
\newline
If $v_{2}\left( \allowbreak x\allowbreak \right) =1$, then $v_{2}\left(
x\left( 2_{E}(\allowbreak x,\allowbreak y)\right) \right) =2$, $v_{2}\left(
x\left( 4(\allowbreak x,\allowbreak y)\right) \right) =0$, $v_{2}\left(
x\left( 8(\allowbreak x,\allowbreak y)\right) \right) =-2$, and $
8(\allowbreak x,\allowbreak y)$ is not an integral point.
\newline
If $v_{2}\left( \allowbreak x\allowbreak \right) \geq 2$, then $v_{2}\left(
x\left( 2_{E}(\allowbreak x,\allowbreak y)\right) \right) =0$, $v_{2}\left(
x\left( 4(\allowbreak x,\allowbreak y)\right) \right) =-2$, and $
4(\allowbreak x,\allowbreak y)$ is not an integral point.
\newline
In all three cases, $8(\allowbreak x,\allowbreak y)$ is not an integral point. 
\newline\newline
Proof of Lemma 2
\newline
We have $k,l\in \mathbb{Z}
$, $k,l>0$, $k\neq l$, $\gcd \left( k,l\right) =1$, and $k$ and $l$ are
square-free.
\newline
We may assume without loss of generality that $p|k$.  Then $v_{p}\left( k\right) =1$ and $v_{p}\left( l\right) =0$.  Since $v_{p}\left( k\right) =1$, $\exists j\in 
\mathbb{Z}
-\left\{ 0\right\} $ such that $k=pj$, $v_{p}\left( j\right) =v_{p}\left(
l\right) =0$.
\newline
Then $a=-3k^{2}l^{2} = -3p^{2}j^{2}l^{2}$ and $b=k^{2}l^{2}\left( k^{2}+l^{2}\right)=p^{2}j^{2}l^{2}\left( p^{2}j^{2}+l^{2}\right)$.
\newline
Also, there exist unique $r,u\in \mathbb{Z}$ such that $x=pu+r$ and $r\in \left\{ 0,1,...p-1\right\} $.  Substituting these values for $a$, $b$, and $x$ into the formula for $x\left( 3\left( x,y\right) \right) $ given by part (ii) of the Proposition, we obtain  $x\left( 3\left( x,y\right) \right) =\frac{r\alpha +p^{6}\left( 64j^{6}l^{12}+p\beta \right) \allowbreak }{
9r\gamma +9p^{6}\delta ^{2}\allowbreak }=\frac{p\epsilon +r^{9}\allowbreak }{
9r\gamma +9p^{6}\delta ^{2}\allowbreak }$, where $\alpha$, $\beta$, $\gamma$, $\delta$, and $\epsilon$ are integers.
If $r=0$, then $x\left( 3\left( x,y\right) \right) =\frac{p^{6}\left(
64j^{6}l^{12}+p\beta \right) \allowbreak }{9p^{6}\delta ^{2}\allowbreak }=
\frac{p\beta +64j^{6}l^{12}}{9\delta ^{2}}$.
If $r\neq 0$, then $x\left( 3\left( x,y\right) \right) =\frac{p\epsilon
+r^{9}\allowbreak }{9r\gamma +9p^{6}\delta ^{2}\allowbreak }$.
In both cases, $p$ does not divide the numerator of $x\left( 3\left(
x,y\right) \right) $ and
$v_{p}\left( x\left( 3\left( x,y\right) \right) \right) \leq \left\{ 
\begin{array}{lll}
-2\text{ } & \text{if } & p=3, \\ 
0 & \text{if } & \text{otherwise.}
\end{array}\right. $
\newline\newline\newline
Proof of Lemma 3
\newline
We have $k,l\in \mathbb{Z}
$, $k,l>0$, $k\neq l$, $\gcd \left( k,l\right) =1$, and $k$ and $l$ are
square-free.  We may assume without loss of generality that $p|k$.  Then $v_{p}\left( k\right) =1$ and $v_{p}\left( l\right) =0$.  By Lemma 2, $v_{p}\left( x\left( 3\left( x,y\right) \right) \right) \leq 0$.  Hence, $\exists w,z\in \mathbb{Z}-\left\{ 0\right\} $ such that $x\left( 3(\allowbreak x,\allowbreak
y)\right) =\frac{w}{z}$ and $v_{p}\left( w\right) =0$.  Since $k^{2}l^{2}$ divides both $a$ and $b$, part (vii) of the Proposition implies that $z^{p^{2}-3}\left( \psi _{p}^{2}-p^{2}\left( \frac{w}{z}\right)
^{p^{2}-1}\right) $ and $z^{p^{2}-2}\left( \phi _{p}-\left( \frac{w}{z}
\right) ^{p^{2}}\right) $ are integers, each of which is divisible by $
k^{2}l^{2}$.  Hence, $x\left( 3p(\allowbreak x,\allowbreak y)\right) =x\left( p\left(
3(\allowbreak x,\allowbreak y)\right) \right) =x\left( p\left( x\left( 3(\allowbreak x,\allowbreak y)\right) ,y\left(
3(\allowbreak x,\allowbreak y)\right) \right) \right)
=x\left( p\left( \frac{w}{z},y\left( 3(\allowbreak x,\allowbreak y)\right)
\right) \right) =\frac{\phi _{p}}{\allowbreak \psi _{p}^{2}}=\frac{
z^{p^{2}-2}\phi _{p}}{z^{p^{2}-2}\allowbreak \psi _{p}^{2}}=\frac{
z^{p^{2}-2}\phi _{p}}{zz^{p^{2}-3}\allowbreak \psi _{p}^{2}}=\frac{
z^{p^{2}-2}\left( \phi _{p}-\left( \frac{w}{z}\right) ^{p^{2}}+\frac{
w^{p^{2}}}{z^{p^{2}}}\right) }{zz^{p^{2}-3}\left( \allowbreak \psi
_{p}^{2}-p^{2}\left( \frac{w}{z}\right) ^{p^{2}-1}+\frac{p^{2}w^{p^{2}-1}}{
z^{p^{2}-1}}\right) }=\frac{z^{p^{2}-2}\left( \phi _{p}-\left( \frac{w}{z}\right)
^{p^{2}}\right) +z^{p^{2}-2}\left( \frac{w^{p^{2}}}{z^{p^{2}}}\right) }{
z\left( z^{p^{2}-3}\left( \allowbreak \psi _{p}^{2}-p^{2}\left( \frac{w}{z}
\right) ^{p^{2}-1}\right) +z^{p^{2}-3}\left( \frac{p^{2}w^{p^{2}-1}}{
z^{p^{2}-1}}\right) \right) }=\frac{k^{2}l^{2}\sigma +\frac{w^{p^{2}}}{z^{2}}
}{z\left( k^{2}l^{2}\tau +\frac{p^{2}w^{p^{2}-1}}{z^{2}}\right) }$, where $\sigma ,\tau \in \mathbb{Z}$.  Since $v_{p}\left( k\right) =1$, $\exists j\in \mathbb{Z}$ with $j\geq 1$ such that $k=pj$ and $v_{p}\left( j\right) =0$.  Then $x\left( 3p(\allowbreak x,\allowbreak y)\right) =\frac{\left( pj\right)
^{2}l^{2}\sigma +\frac{w^{p^{2}}}{z^{2}}}{z\left( \left( pj\right)
^{2}l^{2}\tau +\frac{p^{2}w^{p^{2}-1}}{z^{2}}\right) }=\frac{
p^{2}j^{2}l^{2}\sigma z^{2}+w^{p^{2}}}{p^{2}z\left( j^{2}l^{2}\tau
z^{2}+w^{p^{2}-1}\right) }$ and $v_{p}\left( x\left( 3p\left( x,y\right) \right) \right) \leq -2$.

\end{document}